\documentclass{article}
\usepackage{amsmath, amssymb}
\usepackage{amssymb,pb-diagram,pst-all}
\usepackage{amscd}
\usepackage[dvipdfmx]{graphicx}
\usepackage{tikz}
\usepackage[subrefformat=parens]{subcaption}
\usepackage{fancybox}

\newcommand{\PP}{\mathbb P}

\newcommand{\CC}{\mathbb C}

\newcommand{\mcB}{\mathcal B}

\newcommand{\mcL}{\mathcal L}

\newcommand{\mcQ}{\mathcal Q}

\newcommand{\tor}{{\mathrm {tor}}}

\newcommand{\MW}{\mathop {\rm MW}\nolimits}

\newcommand{\Sing}{\mathop {\rm Sing}\nolimits}

\newtheorem{thm}{Theorem}[section]
\newtheorem{cor}{Corollary}[section]
\newtheorem{prop}{Proposition}[section]
\newtheorem{lem}{Lemma}[section]
\newtheorem{defin}{Definition}[section]
\newtheorem{exmple}{Example}[section]
\newtheorem{rem}{Remark}[section]
\newtheorem{qz}{Question}[section]

\newcommand{\I}{\mathop {\rm I}\nolimits}

\newcommand{\III}{\mathop {\rm III}\nolimits}

\newenvironment{remark}{\begin{rem}\rm }{\end{rem}}

\newcommand{\qed}{\hfill $\Box$}

\newcommand{\proof}{\noindent{\textsl {Proof}.}\hskip 3pt}
\newcommand{\proofend}{\qed \par\smallskip\noindent}

\renewcommand{\thesubparagraph}{\theparagraph.\@arabic\c@subparagraph}
\begin{document}
  
  \begin{center}

{\Large \bf A note on the topology of arrangements for  a smooth plane quartic 
and 
its bitangent lines}

\bigskip

{\bf 
Shinzo BANNAI, Hiro-o TOKUNAGA and  Momoko YAMAMOTO
}

\end{center}
\normalsize

\begin{abstract}
In this paper, we give a Zariski triple of the arrangements for a smooth quartic and its four bitangents. A key criterion to distinguish
the topology of such curves is given by a matrix related to the height pairing of rational points arising
from three bitangent lines.


\medskip


\end{abstract}

\section*{Introduction}

Let $(\mcB_1, \mcB_2)$ be a pair of reduced plane curves.
The pair $(\mcB^1, \mcB^2)$ is said to be a {\it Zariski pair}
 if it satisfies the two conditions as follows:

\begin{enumerate}

\item[(i)] For each $i$, there exists a tubular neighborhood $T(\mcB^i)$ of $\mcB^i$ such
that $(T(\mcB^1), \mcB^1)$ is homeomorphic to $(T(\mcB^2), \mcB^2)$.

\item[(ii)] There exists no homeomorphism between 
$(\PP^2, \mcB^1)$ and $(\PP^2, \mcB^2)$.
\end{enumerate}

An $N$-ple $(\mcB_1, \ldots, \mcB_N)$ is called a Zariski $N$-ple if $(\mcB_i, \mcB_j)$ is a Zariski pair for any $1 \le i < j \le N$.
The first condition for a Zariski pair is replaced by {\it the combinatorics} (or {the combinatorial type}) of $\mcB^i$. For the precise definition
of the combinatorics,  see \cite{act} (It can also be found in \cite{tokunaga14}).  Since the combinatorics is 
more tractable, we always consider the combinatorics rather than the homeomorphism type
 of $\mcB^i$. In \cite{zariski29}, Zariski first found that the topology of
a pair $(\PP^2, \mcB)$  is not determined by the combinatorics of $\mcB$ in the case where $\mcB$ is an irreducible sextic 
with $6$ cusps as its singularities. We refer to \cite{act} for results on Zariski pairs before 2008.  Within these several years,
new approaches to study Zariski pairs for reducible plane curves have been introduced, such as
(a) linking sets (\cite{benoit-jb}), (b) splitting types (\cite{bannai16}), (c) splitting and connected numbers (\cite{shirane16, shirane17}) and (d) the set of subarrangements of $\mcB$ (\cite{bgst, bannai-tokunaga}).

 In \cite{benoit-jb, bgst}, 
Zariski pairs for a smooth cubic and its $k$-inflectional tangents ($k \ge 4$)  were investigated based on the method (d) as above.
This generalizes  E.~Artal's  Zariski pair for a smooth cubic and its three inflectional tangents given in \cite{artal94}.  In \cite{shirane17},  Shirane introduced {\it connected numbers} and generalized  E.~Artal's example to smooth curves of higher degree.

Also, in constructing plane curves which can be candidates for Zariski pairs,  the first and the second authors introduced a new method by
using the geometry of sections and multi-sections of an elliptic surface (\cite{bannai-tokunaga, bannai-tokunaga17, tokunaga14}).  In \cite{bannai16,
bannai-tokunaga, bannai-tokunaga17}, with the methods (b) and (d), they gave some examples for Zariski $N$-plet for arrangements
of curves with low degrees.

In this article, we consider Zariski pairs for a smooth quartic and its bitangents, which can be considered as not only a continuation of
previous studies (e.g., \cite{bgst}), but also as a new point of view for such a classically well-studied object.

A smooth quartic $\mcQ$ and its 28 bitangents have been studied intensively by various authors and there are a lot of results on them. A detailed account of the history of the study of quartic curves and their bitangents can be found in \cite[Chapter 6]{dolg}. 
As for Zariski pairs, however, there do not seem to be any results except a Zariski pair for a smooth quartic and its three bitangents by
E.~Artal and J.~Vall\'es, about which the authors were informed via private communication. In this article, we study such objects
through the Mordell-Weil lattices, the connected numbers and the set of subarrangements.  Here are our main results:

\begin{thm}\label{thm:main}{Consider two combinatorial types of arrangements consisting of a smooth quartic $\mcQ$ and some of its bitangents as follows:

\begin{enumerate}

 \item[(a)] the quartic $\mcQ$  and  three of its bitangent lines which are non-concurrent, 
 
 \item[(b)] the quartic $\mcQ$ and four of its bitangent lines,  none of three of which are concurrent.

\end{enumerate}

Then the following statements hold:

\begin{enumerate} 

\item[(i)] there exists a Zariski pair  for the arrangement (a).

\item[(ii)]  there exists a Zarsiki triple for the arrangement (b)
\end{enumerate}
}
\end{thm}

The first statement  has already been claimed by E.~Artal and J.Vall\'es. Yet we believe that our proof based on 
the theory of Mordell-Weil lattices is different from that of theirs and is new. Hence we believe that it is worthwhile to
present it here.

\begin{remark}
In Corollary \ref{cor:even} we give an upper bound for the number of arrangements, consisting of a smooth quartic and  fixed even number $n$ of bitangent lines, that can be distinguished by our method. Theorem \ref{thm:main} gives an example where the bound is attained in the case where $n=4$. We do not know if the bound is sharp for $n\geq 6$ ($n:$ even).
\end{remark}

In order to explain how we prove Theorem~\ref{thm:main}, We need some preparation.
Let $\mcQ$ be a  smooth quartic
 and choose
 a general point $z_o$ of $\mcQ$. We can associate a rational 
elliptic surface $S_{\mcQ, z_o}$ (see
 \cite[2.2.2]{bannai-tokunaga}, \cite[Section 4]{tokunaga14}) to $\mcQ$ and $z_o$,
 which is given as follows:
 
 \begin{enumerate}
 
 \item[(i)] Let $f_{\mcQ} : S_{\mcQ} \to \PP^2$ be the double cover branched along $\mcQ$. Since $\mcQ$ is smooth,
 $S_{\mcQ}$ is smooth.
 
 
 \item[(ii)] The pencil of lines passing through $z_o$ on $\PP^2$ gives rise to a
  pencil $\Lambda_{z_o}$  of curves of genus $1$ with a unique base 
  point $(f_{\mcQ})^{-1}(z_o)$.
  
  \item[(iii)] Let $\nu_{z_o}: S_{\mcQ, z_o} \to S_{\mcQ}$ be the resolution of the 
  indeterminancy for the rational map induced by $\Lambda_{z_o}$. We denote the induced morphism 
  $\varphi_{\mcQ, z_o} : S_{\mcQ, z_o} \to \PP^1$, which gives a minimal elliptic fibration whose
   generic fiber is denoted by $E_{\mcQ, z_o}$. Note that $E_{\mcQ, z_o}$ is an elliptic curve over $\CC(\PP^1) \cong \CC(t)$.
  The map $\nu_{z_o}$ is a composition of two blowing-ups and the exceptional curve for
  the second blowing-up gives rise to a section $O$ of $\varphi_{\mcQ, z_o}$.
  Note that we have the following diagram:
  \[
\begin{CD}
 S_{\mcQ} @<{\nu_{z_o}}<<S_{\mcQ, z_o} \\
           @VV{f_{\mcQ}}V         @VV{f_{\mcQ, z_o}}V \\
\PP^2 @<< {q_{z_o}}< (\PP^2)_{z_o},
\end{CD}
\]
where $f_{\mcQ, z_o}$ is a double cover induced by the quotient under the involution $[-1]_{\varphi_{\mcQ, z_o}}$ on $S_{\mcQ,z_o}$, which is given by the inversion with respect to the group law 
 on the generic fiber. The morphism $q_{z_o}$ is a composition of 
 two blowing-ups over $z_o$.
  
 \end{enumerate}
 
 If we choose $z_o$ so that the tangent line $l_{z_o}$ at $z_o$ is neither a bitangent line nor a line with intersection multiplicity 4,
 then any bitangent line $L$ of $\mcQ$ gives rise to two section $s_L^{\pm}$. On the generic fiber $E_{S_{\mcQ, zo}}$, we obtain two $\CC(t)$-rational points $\pm P_L$ by
 restricting  these sections to $E_{S_{\mcQ, z_o}}$.
 
 Let us explain how to prove Therorem~\ref{thm:main} (i). Let $L_i$ $(i = 1, 2, 3)$ be three distinct  bitangent lines to $\mcQ$ and let  $\pm P_i$ be the rational
 points obtained from $L_i$ respectively. Put $\triangle = L_1 + L_2 + L_3$.  We then
consider the connected number $c_{f_{\mcQ}}(\triangle)$ (\cite{shirane17} or see \S1) in order to distinguish the topology of $\mcQ+\triangle$.  
In this article, we give a criterion for $c_{f_{\mcQ}}(\triangle)$ to be $1$ or $2$ by using 
a matrix related to the height pairing $\langle P_i, P_j \rangle$ defined by Shioda (\cite{shioda90}).

As for Theorem~\ref{thm:main} (ii), we consider all subarrangements of type $\mcQ+\triangle$ to distinguish the topology of $\mcQ$ and
its four bitangent lines.

The organization of this note is as follows. In  \S 1, we give a brief summary on tools and methods to prove Theorem~\ref{thm:main}.
We give a key criterion in \S 2. Our proof of Theorem~\ref{thm:main} is given in \S 3 where we give an explicit example in the case when
$\mcQ$ is the Klein quartic.

\medskip

{\bf Acknowledgement.} The first author is  partially supported by Grant-in-Aid for Scientific Research C (18K03263).
Also the second author is partially supported by Grant-in-Aid for Scientific Research C (17K05205).
%
%

\section{Preliminaries}

In this section, we introduce various notions which we will use to prove Theorem \ref{thm:main}. The first is the {\it connected number} introduced by T. Shirane in \cite{shirane17}, which will be the key tool in distinguishing the Zariski pair that is claimed to exist in Theorem  \ref{thm:main} (i). Another is the method considered and refined in \cite{bgst}, where the analysis of sub-arrangements effectively distinguishes arrangements with many irreducible components. This method distinguishes the Zariski triple that is claimed to exist in Theorem  \ref{thm:main} (ii). Finally, we introduce the theory of Mordell-Weil lattices, which enables us to conduct the computations needed to apply the above two.

\subsection{Connected Numbers}
%
%
%
%


In \cite{shirane17}, the {\it connected number} is defined for a wide class of varieties, but in this subsection we restate the definition and propositions to fit our setting for the sake of simplicity. The following are simplified versions of \cite[Definition 2.1, Proposition 2.3]{shirane17}.

\begin{defin}
Let  $\phi: X\rightarrow \PP^2$ be a double cover of the projective plane with smooth branch locus $B\subset \PP^2$. Let $C\subset \PP^2$ be a plane curve whose irreducible components are not contained in $B$ such that $C\setminus B$ is connected. Under this setting, the number of connected components of $\phi^{-1}(C\setminus B)$ is called the connected number of $C$ with respect to $\phi$, and will be denoted by $c_{\phi}(C)$.
\end{defin}

Note that we will often omit \lq\lq with respect to $\phi$" when it is apparent from the context. Also, note that since we are considering double covers only, $c_{\phi}(C)=1$ or $2$. The key proposition of connected numbers that will be  used in distinguishing the topology of plane curves is the following:

\begin{prop}
For each $i=1,2$, let $\phi_i: X_i\rightarrow \PP^2$ be a double cover of $\PP^2$ with smooth branch locus $B_i\subset \PP^2$. 
Let $C_1$ be a plane curve whose irreducible components are not contained in $B_1$, such that $C_1\setminus B_1$ is connected. If there exists a homeomorphism $h:\PP^2\rightarrow \PP^2$ with $h(B_1)=B_2$, then $c_{\phi_1}(C_1)=c_{\phi_2}(C_2)$ where $C_2=h(C_1)$.
\end{prop}

\proof
As we are considering double covers only, the assumptions of Proposition 2.3 in \cite{shirane17} are necessarily satisfied if a homeomorphism $h:\PP^2\rightarrow \PP^2$ with $h(B_1)=B_2$ exists. Hence, our statement follows. 
\proofend

\subsection{Distinguishing the embedded topology of plane curves through sub-arrangements}\label{subsec:subarr}

In \cite{bgst},  we formulated a method to study the topology of reducible plane curves via subarrangements. We
here explain its simplified version which fits our case. Let $\mcQ$ be a smooth quartic and $L_i$ ($i = 1, \ldots, 28$)
be its bitangents. Choose a subset $I \subseteq \{1, \ldots, 28\}$ and put $\mcL_I:= \sum_{i\in I} L_i$. Define
\[
\underline{\mathrm{Sub}}_{\triangle}(\mcQ, \mcL_I):=\left . \left \{\mcQ + \sum_{k=1}^3L_{i_k}\right | \forall \{i_1, i_2, i_3\} 
\subset I \right \}.
\]
Define a map $c_I: \underline{\mathrm{Sub}}_{\triangle}(\mcQ, \mcL_I) \to \{1, 2\}$:
\[
c_I : \underline{\mathrm{Sub}}_{\triangle}(\mcQ, \mcL_I) \ni \mcQ + \sum_{k=1}^3L_{i_k} 
\mapsto c_{f_{\mcQ}}\left ( \sum_{k=1}^3L_{i_k}\right ) \in \{1, 2\}
\]
where $f_{\mcQ}$  is the double cover of $\PP^2$ branched along $\mcQ$. Chose two subset $I_1, I_2 \subseteq \{1, \ldots, 28\}$ and let $\mcB_i := \mcQ + \mcL_{I_i}$. If there exists
a homeomorphism $h : (\PP^2, \mcB_1) \to (\PP^2, \mcB_2)$, as $h(\mcQ) = \mcQ$ and $h(\mcL_1) = \mcL_2$
 necessarily hold, it induces a map 
 $h_{\natural} : \underline{\mathrm{Sub}}_{\triangle}(\mcQ, \mcL_{I_1}) \to \underline{\mathrm{Sub}}_{\triangle}(\mcQ, \mcL_{I_2})$ such that
$c_{I_2}=c_{I_1}\circ h_{\natural}$: 
\[
\begin{diagram}
\node{\underline{\mathrm{Sub}}_{\triangle}(\mcQ, \mcL_{I_1})}\arrow{se,t}{c_{I_1}}\arrow{s,l}{h_{\natural}} \\
\node{\underline{\mathrm{Sub}}_{\triangle}(\mcQ, \mcL_{I_2}))}\arrow{e,b}{ c_{I_2}}\node{\{1,2\}} 
\end{diagram}
\]
Hence, likewise \cite[Proposition 1.2]{bgst},  we have the following proposition:

\begin{prop}\label{prop:key}{With the same notation as above, if $\mcB_1$ and $\mcB_2$ have the
same combinatorics and $\sharp c_{I_1}^{-1}(1) \neq \sharp c_{I_2}^{-1}(1)$, then $(\mcB_1, \mcB_2)$ is a
Zariski pair.
}
\end{prop}

\subsection{Elliptic surfaces and Mordell-Weil lattices}

As for basic references about elliptic surfaces and Mordell-Weil lattices, we refer  to \cite{kodaira, miranda-basic, shioda90}. 
In particular, 
for those on rational elliptic surfaces, we refer to \cite{oguiso-shioda}.
In this article, by an {\it elliptic surface}, we always mean the same notion as in \cite{shioda90, bannai-tokunaga}. Namly it means
a smooth projective surface $S$ 
with a relatively minimal genus 1 fibration $\varphi : S \to C$ over a smooth projective curve $C$ with a section 
 $O : C \to S$, which we identify with its image, and at least one singular fiber.
%
%
Let $\Sing(\varphi)=\{ v\in C \mid \text{ $\varphi^{-1}(v)$ is singular }\}$.  For $v \in \Sing(\varphi)$, we put $F_v = \varphi^{-1}(v)$. 
 We denote its irreducible decomposition by 
$ F_v = \Theta_{v, 0} + \sum_{i=1}^{m_v-1} a_{v,i}\Theta_{v,i}$, 
 where $m_v$ is the number of irreducible components of $F_v$ and $\Theta_{v,0}$ is the
unique irreducible component with $\Theta_{v,0}\cdot O = 1$. We call $\Theta_{v,0}$ the {\it identity
 component}.  The classification  of singular fibers is well-known (\cite{kodaira}). We use the Kodaira Notation for the types of singular fibers. 
Let $\MW(S)$ be the set of sections of $\varphi : S \to C$.  We have $\MW(S) \neq \emptyset$ as $O \in \MW(S)$.
By \cite[Theorem 9.1]{kodaira}, $\MW(S)$  is an abelian group with $O$ acting as the zero element. We call
$\MW(S)$  the Mordell-Weil group.   
On the other hand, the generic fiber $E_S$ of
$\varphi : S \to C$ is  a curve of genus $1$ over $\CC(C)$, the rational function field of $C$. The restriction of $O$
to $E$ gives rise to a ${\mathbb C}(C)$-rational point of $E$, and one can regard $E$
as an elliptic curve over ${\mathbb C}(C)$, the restriction of $O$ being the zero element. The group $\MW(S)$ can be identified
with the group of ${\mathbb C}(C)$-rational points $E(\CC(C))$ canonically. 
For $s \in \MW(S)$, we denote the corresponding rational point by $P_s$. 
Conversely,
for an element $P \in E(\CC(C))$,  we denote the corresponding  section by $s_P$.

In \cite{shioda90}, a lattice structure on $(E(\CC(C))/E(\CC(C))_{\tor}$ is defined by using
the intersection pairing on $S$ through $P \mapsto s_P$. 
In particular, $\langle \, , \, \rangle$ denotes the height pairing and $\mbox{\rm Contr}_v$ denotes the contribution term
given in \cite{shioda90} in order to compute $\langle \, , \, \rangle$.

 For the elliptic surface $\varphi_{\mcQ, z_o}: S_{\mcQ, z_o} \to \PP^1$  in the Introduction, $\varphi_{\mcQ, z_o}$ has a unique reducible singular fiber $F_{\infty}$, whose type is
 either $\I_2$ or $\III$ and all other singular fibers are  irreducible. Let $F_{\infty} = \Theta_{\infty, 0} + \Theta_{\infty, 1}$ be
the irreducible decomposition. Then for $P_1, P_2 \in E_{\mcQ, z_o}(\CC(t))$, we have
 \[
 \langle P_1, P_2 \rangle := 1 +s_{P_1}\cdot O + s_{P_2}\cdot O - s_{P_1}\cdot s_{P_2} 
 - \left \{\begin{array}{cc}
           \frac 12 & \mbox{if $\Theta_{\infty, 1}\cdot s_{P_1} = \Theta_{\infty,1}\cdot s_{P_2} = 1$} \\
           0 & \mbox{otherwise}
           \end{array} \right. .           
 \]
 Here, the symbol \lq$\cdot$' denotes the intersection product of divisors.

\section{The height pairing and intersection number of sections}
\subsection{Connected numbers of three bitangents}\label{cnof3b}


Let $\mcQ$ be a smooth plane quartic. We choose homogeneous coordinates $[T, X, Z]$ of $\PP^2$ such that $z_{o}=[0, 1, 0]$ and $Z=0$ is the tangent line of $\mcQ$ at $z_{o}$. Then we may assume that $\mcQ$ is given by a homogeneous polynomial $F_{\mcQ}(T, X, Z)$ of the form
\begin{align*}
F_{\mcQ}(T, X, Z) =ZX^3 + p(T, Z)X^2 + q(T, Z)X+r(T, Z).
\end{align*}
Then the affine part of $\mcQ$, i.e., the part with $Z \neq 0$ is given by
\begin{align*}
F_{\mcQ}(t, x, 1) =x^3 + p(t, 1) x^2 + q(t, 1)x+r(t, 1).
\end{align*}
Then let $\varphi_{\mcQ , z_{o}} : S_{\mcQ, z_{o}} \rightarrow \PP^1$ be the rational elliptic surface as in the Introduction and let $E_{\mcQ , z_{o}}$ be the generic fiber of $\varphi_{\mcQ , z_{o}}$. Then, by \cite[Theorem~10.4]{shioda90}, we have
\begin{align*}
E_{\mcQ , z_{o}}(\CC (t)) \cong E_{7}^{*}, 
\end{align*}
where $E_{7}^{*}$ is the dual lattice of the root lattice $E_{7}$. By \cite{shioda93}, $E_{\mcQ , z_{o}}(\CC (t))$ contains $56$ $\CC(t)$-rational points $P=(x, y)$ of the form:
\begin{align*}
x=at+b, \ y=c t^2 + dt +e.
\end{align*}
Since $-P=(x, -y)$, we denote them by
\begin{align*}
\pm P_{n} = (x_{n}, \pm y_{n}) = \left( a_{n}t+b_{n}, \pm (c_{n}t^2 +d_{n}t +e_{n}) \right) \quad ( n=1, \ldots , 28 ).
\end{align*}
Note that, by \cite[Proposition~4]{shioda93}, $28$ lines $L_n: x_{n}=a_{n}t+b_{n}$ are the $28$ bitangents to $\mcQ$.
As in \S 1, we denote the sections corresponding to $P$ by $s_{P}$. Let $q_{z_o}\circ f_{\mcQ , z_{o}} : S_{\mcQ , z_{o}} \rightarrow \PP^2$ be the map introduced in the Introduction and let $(q_{z_o}\circ f_{\mcQ , z_{o}})^{*}(L_{i})=s_{{i}}^{+}+s_{{i}}^{-}$ $(i = 1,\cdots,28$). Here,  $s_{{i}}^{+}=s_{P_i}$ and $s_{{i}}^{-}=s_{-P_{i}}$. Since $\Theta_{\infty, 1}\cdot s_{\pm P_n} = 1$ and $O\cdot s_{\pm P_n}=0 \ (n= 1, \ldots , 28)$, by the explicit formula for the height pairing, we have the following lemma:

\begin{lem}\label{lem: value of height} For $P_{i}$, $P_{j} \in \{  \pm P_{n} \}$,
\begin{enumerate}
\item[(i)] If $i = j$, then $\langle P_{i} , P_{j} \rangle = \frac{3}{2}$, $s_{P_{i}} \cdot s_{P_{j}}=-1$, and $s_{P_{i}} \cdot s_{-P_{j}}=2$,
\item[(ii)] If $i \neq j$, then 
\begin{enumerate}
\item $\langle P_{i} , P_{j} \rangle =- \frac{1}{2}$ if and only if $s_{P_{i}} \cdot s_{P_{j}}=1$ and $s_{P_{i}} \cdot s_{-P_{j}}=0$,
\item $\langle P_{i} , P_{j} \rangle = \frac{1}{2}$ if and only if $s_{P_{i}} \cdot s_{P_{j}}=0$ and $s_{P_{i}} \cdot s_{-P_{j}}=1$.
\end{enumerate}
\end{enumerate}
\end{lem}


Choose three distinct bitangents $L_{i}$, $L_{j}$, and $L_{k}$ to $\mcQ$. Put $\triangle_{ijk} := L_{i}+L_{j}+L_{k}$. Then, by \S 1, we have connected numbers $c_{f_{\mcQ}}(\triangle_{ijk})=1$ or $2$. 
 From Lemma \ref{lem: value of height}, we classify splitting types of three bitangents via the intersection number of $s_{\pm P_i}$'s. Let the matrix $G(i, j, k)$ be the matrix defined to be two times the Gramm matrix defined by the height pairing of $P_{i}, P_{j}$, and $P_{k}$. The diagonal entries of $G(i,j,k)$ are equal to $3$, and the off-diagonal entries take values $\pm1$. Since $G(i, j, k)$ is a symmetric matrix, then there are $8$ possible choices of $G(i, j, k)$. By the following lemma, the $8$ matrices are classified into two classes depending on $c_{f_{\mcQ}}(\triangle_{ijk})$.
\begin{lem}\label{lem:connected}
\begin{enumerate}
\item[(i)] $c_{f_{\mcQ}}(\triangle_{ijk})=1$ if and only if
\begin{align*}
G(i, j, k) \in 
\left\{ 
\left[
\begin{array}{ccc}
3 & 1 & 1 \\
1 & 3 & 1 \\
1 & 1 & 3
\end{array}
\right]_{,}
\left[
\begin{array}{ccc}
3 & 1 & -1 \\
1 & 3 & -1 \\
-1 & -1 & 3
\end{array}
\right]_{,}
\left[
\begin{array}{ccc}
3 & -1 & 1 \\
-1 & 3 & -1 \\
1 & -1 & 3
\end{array}
\right]_{,}
\left[
\begin{array}{ccc}
3 & -1 & -1 \\
-1 & 3 & 1 \\
-1 & 1 & 3
\end{array}
\right]
 \right\}_{.}
 \end{align*}
\item[(ii)] $c_{f_{\mcQ}}(\triangle_{ijk})=2$ if and only if
\begin{align*}
G(i, j, k) \in 
\left\{ 
\left[
\begin{array}{ccc}
3 & -1 & -1 \\
-1 & 3 & -1 \\
-1 & -1 & 3
\end{array}
\right]_{,}
\left[
\begin{array}{ccc}
3 & -1 & 1 \\
-1 & 3 & 1 \\
1 & 1 & 3
\end{array}
\right]_{,}
\left[
\begin{array}{ccc}
3 & 1 & -1 \\
1 & 3 & 1 \\
-1 & 1 & 3
\end{array}
\right]_{,}
\left[
\begin{array}{ccc}
3 & 1 & 1 \\
1 & 3 & -1 \\
1 & -1 & 3
\end{array}
\right]
\right\}_{.}
\end{align*}
\end{enumerate}
\end{lem}

\proof
We give a proof when
\begin{align*}
G(i, j, k)=
\left[
\begin{array}{ccc}
3 & 1 & 1 \\
1 & 3 & 1 \\
1 & 1 & 3
\end{array}
\right]
\text{or}
\left[
\begin{array}{ccc}
3 & -1 & -1 \\
-1 & 3 & -1 \\
-1 & -1 & 3
\end{array}
\right]_{,}
\end{align*}
since the proof for the other 6 matrices can be done in the same manner.
\begin{enumerate}
\item[(i)] If
\begin{align*}
G(i, j, k)=
\left[
\begin{array}{ccc}
3 & 1 & 1 \\
1 & 3 & 1 \\
1 & 1 & 3
\end{array}
\right]_{,}
\end{align*}
i.e., $2\langle P_{i}, P_{j} \rangle = 2\langle P_{j}, P_{k} \rangle = 2\langle P_{k}, P_{i} \rangle =1$, we have $s_{P_{i}} \cdot s_{P_{j}}=s_{P_{j}} \cdot s_{P_{k}}=s_{P_{k}} \cdot s_{P_{i}}=0$ and $s_{P_{i}} \cdot s_{-P_{j}}=s_{P_{j}} \cdot s_{-P_{k}}=s_{P_{k}} \cdot s_{-P_{i}}=1$ from Lemma \ref{lem: value of height}. Since $\langle \ , \ \rangle$ is symmetric, we obtain $s_{P_{i}} \cdot s_{-P_{j}}=s_{-P_{j}} \cdot s_{P_{k}}=s_{P_{k}} \cdot s_{-P_{i}}=s_{-P_{i}} \cdot s_{P_{j}}=s_{P_{j}} \cdot s_{-P_{k}}=s_{-P_{k}} \cdot s_{P_{i}}=1$. This means that $c_{f_{\mcQ}}(\triangle_{ijk})=1$. 

\item[(ii)]
If 
\begin{align*}
G(i, j, k)=
\left[
\begin{array}{ccc}
3 & -1 & -1 \\
-1 & 3 & -1 \\
-1 & -1 & 3
\end{array}
\right]_{,}
\end{align*}
i.e., $2\langle P_{i}, P_{j} \rangle = 2\langle P_{j}, P_{k} \rangle = 2\langle P_{k}, P_{i} \rangle =-1$, we have $s_{P_{i}} \cdot s_{P_{j}}=s_{P_{j}} \cdot s_{P_{k}}=s_{P_{k}} \cdot s_{P_{i}}=1$ and $s_{P_{i}} \cdot s_{-P_{j}}=s_{P_{j}} \cdot s_{-P_{k}}=s_{P_{k}} \cdot s_{-P_{i}}=0$ from Lemma \ref{lem: value of height}. Hence we obtain $s_{P_{i}} \cdot s_{P_{j}}=s_{P_{j}} \cdot s_{P_{k}}=s_{P_{k}} \cdot s_{P_{i}}=1$ and $s_{-P_{i}} \cdot s_{-P_{j}}=s_{-P_{j}} \cdot s_{-P_{k}}=s_{-P_{k}} \cdot s_{-P_{i}}=1$, i.e., $c_{f_{\mcQ}}(\triangle_{ijk})=2$. 
\end{enumerate}
\proofend

The figures below explain configurations of $(q_{z_o}\circ f_{\mcQ, z_o})^{-1}(\Delta_{ijk}\setminus \mcQ)$ in case (i), (ii) from the proof of Lemma \ref{lem:connected}. Note that the preimages of points on $\mcQ$ are
ignored.
\begin{figure}[h]
\begin{minipage}{0.5\linewidth}
\centering
\scalebox{0.8}{
{\unitlength 0.1in%
\begin{picture}( 27.1400, 24.4600)( -1.4500,-24.5100)%
%
\special{pn 20}%
\special{pa 2300 520}%
\special{pa 2300 2120}%
\special{fp}%
\special{pa 2300 2120}%
\special{pa 2300 2120}%
\special{fp}%
%
\special{pn 20}%
\special{pa 1160 145}%
\special{pa 2548 941}%
\special{fp}%
\special{pa 2548 941}%
\special{pa 2548 941}%
\special{fp}%
%
\special{pn 20}%
\special{pa 560 520}%
\special{pa 560 2120}%
\special{fp}%
\special{pa 560 2120}%
\special{pa 560 2120}%
\special{fp}%
%
\special{pn 20}%
\special{pa 300 1655}%
\special{pa 1688 2451}%
\special{fp}%
\special{pa 1688 2451}%
\special{pa 1688 2451}%
\special{fp}%
%
\special{pn 20}%
\special{pa 2569 1656}%
\special{pa 1180 2450}%
\special{fp}%
\special{pa 1180 2450}%
\special{pa 1180 2450}%
\special{fp}%
%
\special{pn 20}%
\special{pa 1699 156}%
\special{pa 310 950}%
\special{fp}%
\special{pa 310 950}%
\special{pa 310 950}%
\special{fp}%
\put(1.4000,-9.7000){\makebox(0,0){$s_{P_{i}}$}}%
\put(11.2000,-0.7000){\makebox(0,0){$s_{-P_{j}}$}}%
\put(23.0000,-4.3000){\makebox(0,0){$s_{P_{k}}$}}%
\put(27.3000,-15.4000){\makebox(0,0){$s_{-P_{i}}$}}%
\put(17.2000,-25.2000){\makebox(0,0)[lb]{$s_{P_{j}}$}}%
\put(5.6000,-22.0000){\makebox(0,0){$s_{-P_{k}}$}}%
\end{picture}}%
}
\label{fig:connum1}
\subcaption{Case (i), where $c_{f_{\mcQ}}(\triangle_{ijk})=1$.}
\end{minipage}
\begin{minipage}{0.5\linewidth}
\centering
\scalebox{0.8}{
{\unitlength 0.1in%
\begin{picture}( 23.9500, 26.9500)(  7.2500,-31.3500)%
%
\special{pn 20}%
\special{pa 1965 553}%
\special{pa 3022 2380}%
\special{fp}%
\special{pa 3022 2380}%
\special{pa 3022 2380}%
\special{fp}%
%
\special{pn 20}%
\special{pa 1110 2360}%
\special{pa 2169 534}%
\special{fp}%
\special{pa 2169 534}%
\special{pa 2169 534}%
\special{fp}%
%
\special{pn 20}%
\special{pa 1010 2910}%
\special{pa 3120 2910}%
\special{fp}%
\special{pa 3120 2910}%
\special{pa 3120 2910}%
\special{fp}%
%
\special{pn 20}%
\special{pa 1010 2190}%
\special{pa 3120 2190}%
\special{fp}%
\special{pa 3120 2190}%
\special{pa 3120 2190}%
\special{fp}%
%
\special{pn 20}%
\special{pa 1110 3080}%
\special{pa 1580 2270}%
\special{fp}%
\special{pa 1580 2270}%
\special{pa 1580 2270}%
\special{fp}%
%
\special{pn 20}%
\special{pa 1680 2110}%
\special{pa 2180 1248}%
\special{fp}%
\special{pa 2180 1248}%
\special{pa 2180 1248}%
\special{fp}%
%
\special{pn 20}%
\special{pa 2459 2109}%
\special{pa 1958 1248}%
\special{fp}%
\special{pa 1958 1248}%
\special{pa 1958 1248}%
\special{fp}%
%
\special{pn 20}%
\special{pa 3020 3090}%
\special{pa 2554 2278}%
\special{fp}%
\special{pa 2554 2278}%
\special{pa 2554 2278}%
\special{fp}%
\put(10.1000,-23.9000){\makebox(0,0){$s_{P_{i}}$}}%
\put(17.5000,-5.7000){\makebox(0,0)[lb]{$s_{P_{j}}$}}%
\put(32.8000,-21.9000){\makebox(0,0){$s_{P_{k}}$}}%
\put(11.2000,-32.0000){\makebox(0,0){$s_{-P_{i}}$}}%
\put(31.2000,-32.0000){\makebox(0,0){$s_{-P_{j}}$}}%
\put(33.1000,-29.1000){\makebox(0,0){$s_{-P_{k}}$}}%
\end{picture}}%
}
\label{fig:connum2}
\subcaption{Case (ii),  where $c_{f_{\mcQ}}(\triangle_{ijk})=2$.}
\end{minipage}
\end{figure}

By observing the matrices in the two classes above, we have the following Lemmas:
\begin{lem}\label{lem:cnof3b}
Let $m_{ijk}$ be the number of upper-half entries of $G(i,j,k)$ taking values equal to $-1$. Under the above setting,
\begin{enumerate}
\item[(i)] $c_{f_{\mcQ}}(\triangle_{ijk}) = 1$ if and only if $m_{ijk}$  is even,
\item[(ii)] $c_{f_{\mcQ}}(\triangle_{ijk}) = 2$ if and only if  $m_{ijk}$ is odd.\end{enumerate}
\end{lem}

We can restate the above Lemma in terms of determinants. Let $I_3$ be the identity matrix of size $3\times 3$.

\begin{lem}\label{lem:cnof3c}
\begin{enumerate}
\item[(i)]  $c_{f_{\mcQ}}(\triangle_{ijk}) = 1$ if and only if $\det(G(i,j,k)-3I_3)=2$,
\item[(ii)] $c_{f_{\mcQ}}(\triangle_{ijk}) = 2$ if and only if $\det(G(i,j,k)-3I_3)=-2$.
\end{enumerate}
\end{lem}
\subsection{The topology of plane quartic and its $n$ bitangents via sub-arrangement}

Let $\mcQ$ be a smooth plane quartic as in \S \ref{cnof3b} and $L_{1}, \ldots , L_{28}$ be $28$ bitangents to $\mcQ$. Choose a subset $I \subset \{ 1, \ldots , 28 \}$ such that $\sharp I = n \ (4 \leq n \leq 28)$ and put $\mcL_{I}:=\sum_{i \in I}L_{i}$. As in \S \ref{cnof3b}, we obtain the $n \times n$ matrix $G_I$ which is defined as twice of the Gramm matrix defined by the height pairing of $P_{i}$'s $(i \in I)$.

In order to consider the embedded topology of $\mcQ +\mcL_{I}$, we use the connected numbers of sub-arrangements $\underline{\mathrm{Sub}}_{\triangle}(\mcQ, \mcL_I)$. Let $c_{I}$ be the map defined in \S \ref{subsec:subarr} and put
\begin{align*}
m_I:=\sharp\left\{\text{ upper-half entries of $G_I$ equal to $-1$ }\right\}  
\end{align*}
Since $\sharp \underline{\mathrm{Sub}}_{\triangle}(\mcQ, \mcL_I) = \binom{n}{3}$, we have
\begin{align*}
\sharp c_{I}^{-1}(2) = \binom{n}{3} - \sharp c_{I}^{-1}(1).
\end{align*}
Hence, there are $\binom{n}{3}+1$ possible pairs 
\begin{align*}
\left( \ \sharp c_{I}^{-1}(1), \ \sharp c_{I}^{-1}(2) \ \right) = \left( 0, \binom{n}{3} \right) , \ \left( 1, \binom{n}{3} - 1 \right), \ldots , \  \left( \binom{n}{3}, 0 \right).
\end{align*} 

By Proposition \ref{prop:key}, it seems that a Zariski $\left(\binom{n}{3}+1\right)$-ple may exist. However, the following proposition shows that this is not true when $n$ is even.
\begin{prop}\label{prp:m(n-2)}
Under the above setting, there exists a non-negative integer $M$ such that
\begin{align*}
m_I \, (n-2) = 2 \, M + \sharp c_{I}^{-1}(2).
\end{align*}
\end{prop}
\proof
We consider the sum
\[
\sum_{ \{ i_{1}, i_{2}, i_{3} \} \subset I} m_{i_{1}i_{2}i_{3}}
\]
where $m_{i_1i_2i_3}$ is defined as in Lemma \ref{lem:cnof3b}. Let us discribe this sum in two ways.
\begin{itemize}
\item[(I)]
Put $\triangle_{i_{1}i_{2}i_{3}}:=L_{i_{1}}+ L_{i_{2}} + L_{i_{3}}$. By Lemma \ref{lem:cnof3b}, we have
\begin{align*}
c_{f_{\mcQ}}(\triangle_{i_{1}i_{2}i_{3}})=1 & \text{ if and only if } m_{i_{1}i_{2}i_{3}}=0, 2, \\
c_{f_{\mcQ}}(\triangle_{i_{1}i_{2}i_{3}})=2 & \text{ if and only if } m_{i_{1}i_{2}i_{3}}=1, 3.
\end{align*}
Put
\begin{align*}
M_{N}:= \sharp \{ \triangle_{i_{1}i_{2}i_{3}} \ | \ m_{i_{1}i_{2}i_{3}}=N \} \ (N=0, 1, 2, 3).
\end{align*}
Then, we have 
\begin{align*}
\sum_{ \{ i_{1}, i_{2}, i_{3} \} \subset I} m_{i_{1}i_{2}i_{3}} &= \sum_{c_{f_{\mcQ}}(\triangle_{i_{1}i_{2}i_{3}})=1} m_{i_{1}i_{2}i_{3}} + \sum_{c_{f_{\mcQ}}(\triangle_{i_{1}i_{2}i_{3}})=2} m_{i_{1}i_{2}i_{3}} \\
&= 0 \cdot M_{0} + 2 \cdot M_{2} + 1 \cdot M_{1} + 3 \cdot M_{3} \\
&= 2 (M_{2} + M_{3}) + M_{1} + M_{3} \\
&= 2 (M_{2} + M_{3}) + \sharp c_{I}^{-1}(2).
\end{align*}
Define $M$ to be $M_{2}+M_{3}$.

\item[(II)] 
Fix an upper-half entry $g_{kl} (\{k,l\}\subset I)$ with value $-1$. Then, by the definitions of the matrices $G_{I}$ and $G({i_{1}, i_{2}, i_{3}})$, we have 
\begin{center}
$G(i_{1}, i_{2}, i_{3})$ contains $g_{kl}$ as  its entry if and only if $k, l \in \{ i_{1}, i_{2}, i_{3} \}$. 
\end{center}
Put $k=i_{1}, l=i_{2}$ without loss of generality. Then we have
\begin{align*}
& \, \sharp \{ G(i_{1}, i_{2}, i_{3}) \ | \ G(i_{1}, i_{2}, i_{3}) \text{ contains } g_{kl} \text{ as its entry } \} \\
= & \, \sharp \{ \{ i_{1}, i_{2}, i_{3} \} \subset I \ | \ k=i_{1}, l=i_{2} \} \\
= & \, n-2.
\end{align*}
Hence, by the definitions of $m_I$ and $m_{i_{1}i_{2}i_{3}}$, we obtain
\begin{align*}
\sum_{\{ i_{1}, i_{2}, i_{3} \} \subset I}m_{i_{1}i_{2}i_{3}}=m_I(n-2).
\end{align*}
\end{itemize}
\proofend

\begin{cor}\label{cor:even}
If $n$ is an even number, then the number of possible pairs of $\left( \, \sharp c_{I}^{-1}(1), \sharp c_{I}^{-1}(2) \, \right)$ is less than or equal to $\frac{1}{2}\binom{n}{3}+1$. 
\end{cor}
\proof
If $n$ is an even number, $\sharp c_{I}^{-1}(2)$ is also even by Proposition \ref{prp:m(n-2)}. Hence, we have
\begin{align*}
( \, \sharp c_{I}^{-1}(1), \sharp c_{I}^{-1}(2) \, ) = \left( 0, \binom{n}{3} \right), \left( 2, \binom{n}{3}-2 \right), \ldots , \left( \binom{n}{3}-2, \, 2 \right), \left( \binom{n}{3}, \, 0 \right).
\end{align*}
\proofend

\begin{rem}{\rm 
We do not know if $\frac{1}{2}\binom{n}{3}+1$ is a strict upper bound for general even $n$. In Section \ref{sec:example}, however,  we give an example attaining $\frac{1}{2}\binom{n}{3}+1$ when $n=4$.} 
\end{rem}

\section{Examples}\label{sec:example}


Let $\mcQ$ be the Klein quartic given by the affine equation:
\begin{align*}
F(t, x):=x^3+t^{3}x+t.
\end{align*}
Then the generic fiber $E_{\mcQ, z_{o}}$ of $\varphi_{\mcQ , z_{o}}$ is $y^2 = F(t, x)$.
From \cite[Section 4]{shioda93}, the $28$ bitangents of $\mcQ$ are given by the following equations:
\begin{align*}
& L_{0, j} :x_{0, j}(t)=-\zeta^{j}t-\zeta^{3j}, 
& L_{1, j}:x_{1, j}(t)=-\zeta^{j}\varepsilon_{1}^2 t-\zeta^{3j}\varepsilon_{3}^{-2}, \\
& L_{2,j}:x_{2, j}(t)=-\zeta^{j}\varepsilon_{2}^2 t-\zeta^{3j}\varepsilon_{1}^{-2}, 
& L_{3,j}:x_{3, j}(t)=-\zeta^{j}\varepsilon_{3}^2 t-\zeta^{3j}\varepsilon_{2}^{-2},
\end{align*}
where $j=0, \ldots , 6, \ \zeta = e^{\frac{2 \pi i}{7}}, \ \varepsilon_{1}=\zeta + \zeta^{-1}, \ \varepsilon_{2}=\zeta^2 + \zeta^{-2}, \ \varepsilon_{3}=\zeta^{4}+\zeta^{-4}$. Put
\begin{align*}
& L_{1}:= L_{0, 0}, \ L_{2}:=L_{1, 0}, \ L_{3}:= L_{1, 1}, \ L_{4}:=L_{3, 3}, \\
& L_{5}:=L_{1, 6}, \ L_{6}:=L_{3, 4}, \ L_{7}:=L_{2, 5}.
\end{align*}
and rational points of $E_{\mcQ, z_{o}}$ defined by $L_{1}, \ldots , L_{7}$;
\begin{align*}
& P_{1}:=\left( \,  x_{0, 0}(t), \ y_{1}(t) \, \right),
P_{2}:=\left( \, x_{1, 0}(t) , \ y_{2}(t) \, \right), 
P_{3}:=\left( \, x_{1, 1}(t), \ y_{3}(t) \, \right),
P_{4}:=\left( \, x_{3, 3}(t), \ y_{4}(t) \, \right), \\
& P_{5}:=\left( \, x_{1, 6}(t), \ y_{5}(t) \, \right),
P_{6}:=\left( \, x_{3, 4}(t), \ y_{6}(t) \, \right),
P_{7}:=\left( \, x_{2, 5}(t), \ y_{7}(t) \, \right).
\end{align*}
Here, 
\begin{align*}
& y_{1}(t)= \sqrt{-1}(t^2 +t+1), \ y_{2}(t)=\sqrt{-1} \, \varepsilon_{1}(t^2 + a_{1}(\zeta) t + b_{1}(\zeta)), \\
& y_{3}(t)= \sqrt{-1}\zeta^{4}\varepsilon_{1}(t^2+\zeta^{2}a_{1}(\zeta) t + \zeta^{4}b_{1}(\zeta)), \ y_{4}(t)= \sqrt{-1}\zeta^{5}\varepsilon_{3}({t}^{2}+ a_{3}(\zeta)t + b_{3}(\zeta) ), \\
& y_{5}(t)=\sqrt{-1}\zeta^{3}\varepsilon_{1}(t^2+ \zeta^{5} a_{1}(\zeta) t+ \zeta^{3}b_{1}(\zeta) ), \ y_{6}(t)=\sqrt{-1}\zeta^{2}\varepsilon_{3}({t}^{2}+\zeta^{2}a_{3}(\zeta)t + \zeta^{4}b_{3}(\zeta)), \\
& y_{7}(t) = \sqrt{-1}\zeta^{6}\varepsilon_{2}({t}^{2}+( \zeta^2 +2 +2 \zeta^6 + \zeta^4 +4 \zeta^3 )t + \zeta^5 + 3 \zeta^3 +3 \zeta^2 + 1 +3 \zeta^6), \\
\intertext{where}
& a_{1}(\zeta) = 2 \zeta^5+\zeta^4+\zeta^3+2 \zeta^2+4, \ b_{1}(\zeta) = 3 \zeta^5+\zeta^4+\zeta^3+3 \zeta^2+3, \\
& a_{3}(\zeta) = 2 \zeta^5 + \zeta^4 + \zeta + 2 + 4 \zeta^6, \ b_{3}(\zeta) = 3 \zeta^4 + \zeta^3 + 1 + 3 \zeta^6 + 3 \zeta^5. 
\end{align*}
Note that for $L_1,\ldots, L_7$, no three lines are concurrent and $\mcQ+\mcL_I$ ($I\subset\{1,\ldots,7\}$) all have the same combinatorics for fixed $\sharp I$.
\begin{itemize}
\item \underline{A Zariski pair for $\mcQ$ and its three bitangents} \\

We put 
\[
\mcB^{1}:= \mcQ + L_{1} + L_{2} + L_{3},\quad
\mcB^{2}:= \mcQ + L_{1} + L_{2} + L_{4}.
\]
For $P_{1}$ and $P_{2}$, consider sections $s_{P_1}$ and $s_{P_2}$ as curves in the affine part of $S_{\mcQ, z_o}$ given by $s_{P_{1}}:=(x_{0, 0}(t), y_{1}(t))$ and $s_{P_{2}}:=(x_{1, 0}(t), y_{2}(t))$ with parameter $t$. Then $x_{0, 0}(t)=x_{1, 0}(t)$ and $y_{1}(t)=y_{2}(t)$ has a unique solution $t=\zeta + \zeta^{-1}$, which implies $s_{P_{1}} \cdot s_{P_{2}}=1$. In the same way, we obtain $s_{P_{2}} \cdot s_{P_{3}}=s_{P_{3}} \cdot s_{P_{1}}=s_{P_{1}} \cdot s_{P_{4}}=1$ and $s_{P_{2}} \cdot s_{P_{4}}=0$. Hence, we have
\[
G(1, 2, 3)=
\left[
\begin{array}{ccc}
3 & -1 & -1 \\
-1 & 3 & -1 \\
-1 & -1 & 3
\end{array}
\right]_{,} \quad
G(1, 2, 4)=
\left[
\begin{array}{ccc}
3 & -1 & -1 \\
-1 & 3 & 1 \\
-1 & 1 & 3
\end{array}
\right]_{.}
\]
By Lemma \ref{lem:cnof3b}, we have $c_{f_{\mcQ}}(\triangle_{123})=2$ and $c_{f_{\mcQ}}(\triangle_{124}) = 1$, then $(\mcB^{1}, \mcB^{2})$ is a Zariski pair.

\item \underline{A Zariski triple for $\mcQ$ and its four bitangents} \\

We set $I_{1}:= \{ 1, 2, 3, 5 \}$, $I_{2}:= \{ 1, 2, 3, 6 \}$, $I_{3}:= \{ 1, 2, 4, 7 \}$ and put
\begin{align*}
\mcB^{k}:= \mcQ + \mcL_{I_{k}} \ (k = 1, 2, 3).
\end{align*}
As above, we have $G(1, 2, 3) = G(1, 2, 5) = G(1, 3, 5)$, $G(1, 2, 4)=G(1, 2, 7)=G(1, 4, 7)$ and 
\begin{align*}
& G(1, 2, 6) = G(1, 3, 6) =
\left[
\begin{array}{ccc}
3 & -1 & -1 \\
-1 & 3 & 1 \\
-1 & 1 & 3
\end{array}
\right]_{,} \\
& G(2, 3, 6)=
\left[
\begin{array}{ccc}
3 & -1 & 1 \\
-1 & 3 & 1 \\
1 & 1 & 3
\end{array}
\right]_{,} \quad
G(2, 4, 7)=
\left[
\begin{array}{ccc}
3 & 1 & 1 \\
1 & 3 & 1 \\
1 & 1 & 3
\end{array}
\right]_{.}
\end{align*}
Hence, we obtain
\[
\left( \, \sharp c_{I_{1}}^{-1}(1), \sharp c_{I_{1}}^{-1}(2) \, \right) = ( 4, 0 ),\,
\left( \, \sharp c_{I_{2}}^{-1}(1), \sharp c_{I_{2}}^{-1}(2) \, \right) = ( 2, 2 ), \,
\left( \, \sharp c_{I_{3}}^{-1}(1), \sharp c_{I_{3}}^{-1}(2) \, \right) = ( 0, 4 ).
\]
By Proposition \ref{prop:key}, $(\mcB^{1}, \mcB^{2}, \mcB^{3})$ is a Zariski triple.
\end{itemize}

The existence of the above examples gives a proof to Theorem \ref{thm:main}.


\noindent Shinzo BANNAI\\
National Institute of Technology, Ibaraki College, 866 Nakane, Hitachinaka-shi, Ibaraki-Ken 312-8508 JAPAN, 
{\tt sbannai@ge.ibaraki-ct.ac.jp}\\

\noindent Hiro-o TOKUNAGA, Momoko YAMAMOTO\\
Department of Mathematics Sciences\\
Tokyo Metropolitan University, 1-1 Minami-Ohsawa, Hachiohji 192-0397 JAPAN,\\ 
{\tt tokunaga@tmu.ac.jp}, {\tt yamamoto-momoko@ed.tmu.ac.jp}
%
%

{\tt }

\vspace{0.5cm}
      
 \end{document}